\documentclass[12pt]{article}
\usepackage{amsmath,amsfonts,amssymb}

\usepackage{geometry}
\usepackage{amsmath}
\usepackage{amsfonts}
\usepackage{amssymb}

 \newtheorem{theorem}{Theorem}

\author{Sergey V. Galaev}

\title{Almost contact metric structures defined by an $N$-prolonged connection}

\begin{document}

\maketitle

\begin{abstract} On a manifold with an almost contact metric
structure $(\varphi,\vec\xi,\eta,g,X,D)$ the notions of the
interior and the $N$-prolonged connections are introduced. Using
the $N$-prolonged connection, a new almost contact metric
structure is defined on the distribution $D$. The properties of
this structure are studied.


\end{abstract}

\section[1]{Introduction}
 The study of the geometry of the tangent  bundles was initiated
 in the fundamental paper [1] by Sasaki published in 1958.
 Using a Riemannian metric $g$ on a manifold $X$ Sasaki defines a Riemannian metric $G$
 on the tangent bundle $TX$ of the manifold $X$. This construction is grounded on the natural splitting
(that takes a place due to the existence of the Levi-Civita
connection) of the tangent bundle $TTX$ of the manifold $TX$ into
the direct sum  of the vertical and horizontal distributions, the
fibers of these distributions are isomorphic to the fibers of the
distribution $TX$. The odd analogue of the tangent bundle is a
distribution $D$ of an almost contact metric structure
$(\varphi,\vec{\xi},\eta,g)$. Similarly to the bundle $TTX$, the
bundle $TD$ due to a connection over the distribution [2] (and
later an $N$-prolonged connection, i.e. connection in the vector
bundle $(X,D)$) splits into the direct sum of the vertical and
horizontal distributions. In [2,3] it is shown that on the
manifold $D$ can be defined in a natural way an almost contact
metric structure allowing e.g. to give an invariant character to
the analytical description of the mechanics with constraints. In
[3], on the manifold $D$ the geodesic pulverization of the
connection over the distribution is defined; this is an analogue
of the geodesic pulverization, defined on the space of the tangent
bundle $TX$ and having the following physical interpretation: the
projections of the integral curves of the geodesic pulverization
of the connection over the distribution coincide with the
admissible geodesics  (the trajectories of the mechanical system
with constrains).

The present work is an introduction to the geometry of prolonged
almost contact metric structures and it is dedicated to the
development of the following two ideas: the idea of the
generalization  of the Sasaki construction [1] for the case of odd
dimension, and the idea of the extension of the  interior
connection.

The paper has the following structure. The second section consists
of three subsections;  the first of them contains short
information about the interior geometry of almost contact metric
spaces. This stuff can be found in [4].

In section 2.2, the notion of the $N$-prolonged metric connection
is introduced. The interior connection defines the parallel
transport of the admissible vectors (i.e. the vectors belonging to
the distribution $D$) along admissible curves. Each corresponding
$N$-prolonged connection is a connection in the vector bundle $(D,
\pi,X)$ defined by the  interior  connection and an endomorphism
$N:D\to D$. The choice of the endomorphism $N:D\to D$ effects the
properties of the prolonged  connection and also the properties of
the (prolonged) almost contact metric structure appearing on the
total space $D$ of the vector bundle $(D, \pi , X)$. The central
of this subsection is the theorem about the existence and
uniqueness of the $N$-prolonged metric  connection with zero
torsion. In section 2.3, the relation of the interior and
prolonged connections with the known connections  appearing on
almost contact metric spaces is shown.

In the third section, on the manifold $D$ with an prolonged metric
connection, the prolonged almost contact metric structure is
defined. In subsection 3.1, the properties of the prolonged almost
contact metric structure are investigated. In subsection  3.2, a
negative answer to the question about the possibility of the
isometrical imbedding of the manifold $D$ to the manifold $TX$
with the Sasaki metric is given.


\section[2]{Interior and $N$-prolonged connections}

\subsection[2.1]{Preliminaries on the interior geometry of almost
contact metric spaces}

Let $X$  be a smooth manifold of  an odd dimension $n$. Denote by
$\Gamma TX$ the $C^{\infty}(X)$-module of smooth vector fields on
$X$. All manifolds, tensors and other geometric objects will be
assumed to be smooth of the class $C^{\infty}$.  An almost contact
metric structure on
 $X$ is an aggregate
 $(\varphi, \vec{\xi}, \eta, g)$ of  tensor fields on $X$, where $\varphi$ is a tensor field of type
 $(1, 1)$, which is called the structure
 endomorphism, $\vec{\xi}$ and $\eta$  are a vector and a covector,
 which are called the structure vector and the contact form,
 respectively, and  $g$ is a (pseudo-)Riemannian metric.
 Moreover,
$$\eta(\vec{\xi})=1,\quad \varphi(\vec{\xi})=0,\quad \eta \circ
\varphi=0,$$
$$\varphi^2\vec{X}=-\vec{X}+\eta(\vec{X})\vec{\xi},\quad
g(\varphi\vec{X},\varphi\vec{Y})=g(\vec{X},\vec{Y})-\eta(\vec{X})\eta(\vec{Y}),$$
$$d\eta(\vec\xi,\vec X)=0$$

for all $\vec{X}, \vec{Y} \in \Gamma TM$. The skew-symmetric
tensor
 $\Omega(\vec{X}, \vec{Y})=g(\vec{X},
\varphi\vec{Y})$ is  called the fundamental form of the structure.
A manifold with a fixed almost contact metric structure is called
an almost contact metric manifold. If
 $\Omega=d\eta$ holds, then the almost contact metric structure is called a contact metric structure.
An almost contact metric structure is called normal if
$$N_{\varphi}+2d\eta\otimes\vec{\xi}=0,$$ where $N_{\varphi}$ is the
Nijenhuis torsion defined for the tensor $\varphi$. A normal
contact metric structure is called a Sasakian structure. A
manifold with a given Sasakian structure is called a Sasakian
manifold. Let $D$ be the smooth distribution of codimension 1
defined by the form  $\eta$, and $D^\bot={\rm span}(\vec{\xi})$ be
the closing of $D$. If  the restriction of the 2-form
$\omega=d\eta$ to the distribution $D$  is non-degenerate, then
the vector $\vec{\xi}$ is uniquely defined by the condition
$$\eta(\vec{\xi})=1,\quad {\rm ker\,} \omega={\rm span}
(\vec{\xi}),$$ and  it is called the Reeb vector field. In this
subsection we will pay more attention to the so called almost
contact K\"ahlerian spaces [5], other basic classes of of almost
contact metric spaces will be considered in the next section.

An almost contact metric structure is called almost normal, if it
holds
\begin{equation}\label{eq1}N_{\varphi}+2(d\eta\circ\varphi)\otimes\vec{\xi}=0.\end{equation}
In what follows, an almost normal almost contact metric space will
be called {\it an almost contact K\"ahlerian space}  if its
fundamental form is closed. An almost contact metric space is
called {\it almost K-contact metric space} if $L_{\vec\xi}g=0$.
The last equality is usually used in the case, when the form
$\omega$ has the maximal rank, then the corresponding space is
called K-contact.

An almost normal contact metric structure is obviously a Sasakian
structure. Sasakian manifolds are popular among the researchers of
almost contact metric spaces by the following two reasons. On one
hand, there exist a big number of interesting and deep examples of
Sasakian structures, on the other hand, the Sasakian manifolds
have very important and natural properties. In the same time, the
almost contact K\"ahlerian spaces inherit many important
properties of the Sasakian spaces, this turns out to be very
essential in the cases when an almost contact metric space can not
in principle be a Sasakian space [6].

We say that a coordinate chart $K(x^\alpha)$
$(\alpha,\beta,\gamma=1,...,n,\, a,b,c,e=1,...,n-1)$ on a manifold
$X$ is adapted to the non-holonomic manifold  $D$ if
$$D^{\bot}={\rm span}\left(\frac{\partial}{\partial
x^{n}}\right)$$ holds [4]. Let $$P:TX\rightarrow D$$ be the
projection map defined by the decomposition $$TX=D\oplus
D^{\bot}$$ and let $K(x^{\alpha})$ be an adapted coordinate chart.
Vector fields
$$P(\partial_{a})=\vec{e}_{a}=\partial_{a}-\Gamma^{n}_{a}\partial_{n}$$
are linearly independent, and linearly generate the system $D$
over the domain of the definition of the coordinate map: $$D={\rm
span}(\vec{e}_{a}).$$ Thus we have on $X$ the non-holonomic field
of bases  $$(\vec{e}_{a},\partial_{n})$$ and the corresponding
field of cobases
$$(dx^a,\theta^{n}=dx^{n}+\Gamma^{n}_{a}dx^{a}).$$
 It can be checked directly that
$$[\vec{e}_{a},\vec{e}_{b}]=M^{n}_{ab}\partial_{n},$$ where the
components  $M^{n}_{ab}$ form the so-called tensor of
non-holonomicity [7]. Under assumption that for all adapted
coordinate systems it holds  $\vec{\xi}=\partial_{n}$, the
following equality takes the place
$$[\vec{e}_{a},\vec{e}_{b}]=2\omega_{ba}\partial_{n},$$ where
$\omega=d\eta$.  We say also that the basis
$$\vec{e}_{a}=\partial_{a}-\Gamma^{n}_{a}\partial_{n}$$ is adapted,
as the basis defined by an adapted coordinate map. Note that
$$\partial_n\Gamma^n_a=0$$.
Let $K(x^\alpha)$ and $K(x^{\alpha'})$ be adapted charts, then
under the condition $$\vec\xi=\partial_n$$ we get the next
formulas for the coordinate transformation:
$$x^\alpha=x^{\alpha}(x^{\alpha'}),\quad
x^n=x^{n'}+x^{n}(x^{\alpha'}).$$

A tensor field of type $(p,q)$ defined on an almost contact metric
manifold is called admissible (to the distribution $D$) if in
adapted coordinate map it looks like
$$
t=t^{a_{1},...,a_{p}}_{b_{1},...,b_{q}}\vec{e}_{a_{1}}\otimes...\otimes\vec{e}_{a_{p}}\otimes
dx^{b_{1}}\otimes...\otimes dx^{b_{q}}.
$$
From the definition of an almost contact structure it follows that
the field of endomorphisms $\varphi$ is an admissible tensor field
of type $(1, 1)$. The field of endomorphisms  $\varphi$ we call an
admissible almost complex structure, taking into the account its
properties. The 2-form $\omega=d\eta$ is also an admissible tensor
field and it is natural to call it an admissible symplectic form.

The transformation of the components of an admissible tensor field
in adapted coordinates satisfies the following low:
$$t^a_b=A^a_{a'}A_b^{b'}t_{b'}^{a'},$$
where $A^a_{a'}=\frac{\partial x^a}{\partial x^{a'}}.$

{\bf Remark 1.} From the formulas for the transformation of the
components of an admissible tensor field it follows that the
derivatives $\partial_nt^a_b$ are again components of an
admissible field. Moreover, the vanishing of $\partial_nt^a_b$
does not depend on the choice of adapted coordinates. The last
statement also follows from the equality $(L_{\vec
\xi})^a_b=\partial_nt^a_b$.

{\bf Remark 2.} An admissible tensor structure satisfying
$\partial_nt^a_b=0$ we will call projectible (in the literature
there are other names for the structures with similar properties:
basic, semi basic, etc.). In what follows we will see that
admissible projectible structures can be naturally considered as
structures defined on a submanifold of a smaller dimension.

Using adapted coordinates we introduce the following admissible
tensor fields: $$h^a_b=\frac{1}{2}\partial_n\varphi^a_b,\quad
C_{ab}=\frac{1}{2}\partial_ng_{ab},\quad C^a_b=g^{da}C_{db},\quad
\psi^b_a=g^{db}\omega_{da}.$$ We denote by $\tilde\nabla$ and
$\tilde\Gamma^\alpha_{\beta\gamma}$ the Levi-Civita connection and
the Christoffel symbols of the metric $g$. The proof of the
following theorem follows from direct computations.

\begin{theorem} The Christoffel symbols of the Levi-Civita connection
of an almost contact metric space with respect to adapted
coordinates are the following:
$$\tilde\Gamma^c_{ab}=\Gamma^c_{ab},\quad \tilde
\Gamma^n_{ab}=\omega_{ba}-C_{ab},\quad \tilde \Gamma^b_{an}=\tilde
\Gamma^b_{na}=C^b_a-\psi^b_a,\quad \tilde \Gamma^n_{na}=\tilde
\Gamma^a_{nn}=0,$$ where $$\Gamma^a_{bc}=\frac{1}{2}g^{ad}(\vec
e_b g_{cd}-\vec e_c g_{bd}-\vec e_d g_{bc}).$$
\end{theorem}

\subsection[2.2]{$N$-prolonged metric connection}

An intrinsic linear connection on a manifold with an almost
contact metric structure [4] is defined as a map
$$\nabla:\Gamma D \times \Gamma D \rightarrow \Gamma D $$ that
satisfies the following conditions:
\begin{align*}
1)\quad &
\nabla_{f_1\vec{u}_1+f_2\vec{u}_2}=f_1\nabla_{\vec{u}_1}+f_2\nabla_{\vec{u}_2};\\
2)\quad &
\nabla_{\vec{u}}f\vec{v}=f\nabla_{\vec{u}}\vec{v}+(\vec{u}f)\vec{v},
\end{align*} where $\Gamma D$ is the module of admissible vector fields. The
Christoffel symbols are defined by the relation $$
\nabla_{\vec{e}_{a}}\vec{e}_{b}=\Gamma^{c}_{ab}\vec{e}_{c}. $$

The torsion $S$ of the intrinsic linear connection is defined by
the formula
$$
S(\vec{X},\vec{Y})=\nabla_{\vec{X}}\vec{Y}-\nabla_{\vec{Y}}\vec{X}-p[\vec{X},\vec{Y}].
$$ Thus with respect to an adapted coordinate system it holds  $$
S^{c}_{ab}=\Gamma^{c}_{ab}-\Gamma^{c}_{ba}. $$

The action of an intrinsic linear connection can be extended in a
natural way to arbitrary admissible tensor fields. An important
example of an intrinsic linear connection is the intrinsic metric
connection that is uniquely defined by the conditions $\nabla g=0$
and $S=0$ [7]. With respect to the adapted coordinates it holds
\begin{equation}\label{eq2} \Gamma^a_{bc}=\frac{1}{2}g^{ad}(\vec e_b g_{cd}-\vec e_c g_{bd}-\vec e_d g_{bc}).
\end{equation}
Note that $\Gamma^a_{bc}=\tilde \Gamma^a_{bc}$ (see Theorem 1).

 In the same way as
a linear connection on a smooth manifold, an intrinsic connection
can be defined by giving a horizontal distribution over the total
space of some vector bundle. In the case of the interior
connection,  the role of such bundle plays the distribution $D$.
 One says that over a distribution  $D$ a connection
is given if the distribution $$\tilde{D}=\pi^{-1}_{*}(D),$$ where
$\pi:D \rightarrow X$ is the natural projection, can be decomposed
into a direct some of the form $$\tilde{D}=HD \oplus VD,$$ where
$VD$ is the vertical distribution on the total space $D$.

Let us introduce a structure of a smooth manifold on $D$. This
structure is defined in the following way. To each adapted
coordinate chart  $K(x^\alpha)$ on the manifold  $X$ we put in
correspondence the coordinate chart
$\tilde{K}(x^{\alpha},x^{n+\alpha})$ on the manifold  $D$, where
$x^{n+\alpha}$ are the coordinates of an admissible vector with
respect to the basis
$$\vec{e}_{a}=\partial_{a}-\Gamma^{n}_{a}\partial_{n}.$$ The
constructed over-coordinate chart will be called adapted. Thus the
assignment of a connection over a distribution is equivalent to
the assignment of the object $$G^{a}_{b}(X^{a},X^{n+a})$$  such
that
$$HD={\rm span}(\vec{\epsilon}_{a}),$$ where
$$\vec{\epsilon}_{a}=\partial_{a}-\Gamma^{n}_{a}\partial_{n}-G^{b}_{a}\partial_{n+b}.$$
If it holds
$$G^{a}_{b}(x^{a},x^{n+a})=\Gamma^{a}_{bc}(x^{a})x^{n+c},$$
then  the connection over the distribution $D$  is defined by  the
intrinsic linear connection. In [2], the notion of the prolonged
connection was introduced. The prolonged
 connection is always considered as a connection over a distribution and it
  is defined by the decomposition
$$TD=\tilde{HD} \oplus VD,$$ where $HD \subset \tilde{HD}$.
The prolonged connection is a connection in a vector bundle. As it
follows from the definition of the prolonged connection, for its
assignment (under the condition that a connection on the
distribution is already defined) it is enough to define a vector
field $\vec u$ on the manifold $D$ that has the following
coordinate form: $\vec u=\partial_n-N^a_bx^{n+b}\partial_{n+a}$,
where the endomorphism $N:D\to D$ can be chosen in an arbitrary
way. We call the torsion of the of the prolonged connection the
torsion of the initial connection. In what follows we call a
prolonged connection an $N$-prolonged connection.

In [7] the admissible tensor field
$$R(\vec u,\vec v)\vec w=\nabla_{\vec u}\nabla_{\vec v}\vec w-\nabla_{\vec v}\nabla_{\vec u}\vec w-
\nabla_{p[\vec u,\vec v]}\vec w$$ is called by Wagner the first
Schouten curvature tensor. With respect to the adapted coordinates
it holds $$R^a_{bcd}=2\vec
e_{[a}\Gamma^d_{b]c}+2\Gamma^d_{[a||e||}\Gamma^e_{b]c}.$$ If the
distribution $D$ does not contain any integrable subdistribution
of dimension $n-2$, then the Schouten curvature tensor  is zero if
and only if the parallel transport of admissible vectors does not
depend on the curve [7]. We say that the Schouten tensor is the
curvature tensor of the distribution $D$. If this tensor is zero,
we say that the distribution $D$ is a zero-curvature distribution.
Note that the partial derivatives
$\partial_n\Gamma^a_{bc}=P^a_{bc}$ are components of an admissible
tensor field [7].

{\bf Remark 3.} In the case of (almost) K-contact spaces the
Schouten tensor has the same properties as the Riemannian tensor
of a manifold. In general this is not true. The vector fields
$$\vec\epsilon_a=\partial_a-\Gamma^n_a\partial_n-\Gamma^b_{ac}x^{n+c}\partial_{n+b},\quad\vec u=\partial_n-N^a_bx^{n+b}\partial_{n+a},\quad\partial_{n+a}$$
define on $D$ a non-holonomic (adapted) field of bases, and the
forms
 $$dx^a,\quad\Theta^n=dx^n+\Gamma^n_adx^a,\quad \Theta^{n+a}=dx^{n+a}+\Gamma^a_{bc}x^{n+c}dx^b+N^a_bx^{n+b}dx^n$$
 define the corresponding field of cobases.
The following structure equations can be obtained:
\begin{align} \label{eq3}
\left[{\vec{\varepsilon }}_a,{\vec{\varepsilon }}_b\right]&=2{\omega }_{ba}\vec u+
x^{n+d}(2\omega_{ba}N^c_d+R^c_{bad}){\partial }_{n+c},\\
\label{eq4}
[{\vec{\varepsilon }}_a,\vec u]&=x^{n+d}(\partial_n \Gamma^c_{ad}-\nabla_aN^c_d){\partial }_{n+c},\\
\label{eq5} \left[{\vec{\varepsilon }}_a,{\partial
}_{n+b}\right]&={\Gamma }^c_{ab}{\partial }_{n+c}.
\end{align}

From \eqref{eq3} and \eqref{eq4} we obtain the the expression for
the curvature tensor of the prolonged connection

\begin{align} \label{eq6} K(\vec u,\vec v)\vec w&=2\omega(\vec
u,\vec v)N\vec w+R(\vec u,\vec v)\vec w,\\
\label{eq7} K(\vec \xi,\vec u)\vec v&=P(\vec u,\vec
v)-(\nabla_{\vec u}N)\vec v,\end{align} where $\vec u,\vec v,\vec
w\in\Gamma D$.

\begin{theorem} There exists an $N$-prolonged metric connection
uniquely determined by the following conditions:

1. $\vec Zg(\vec X,\vec Y)\vec Z=g(\nabla_{\vec Z}\vec X,\vec
Y)+g(\vec X,\nabla_{\vec Z}\vec Y)$ (metricity property),

2. $\nabla_{\vec X}\vec Y - \nabla_{\vec Y}\vec X - p[\vec X,\vec
Y]=0$ (connection is torsion-free),

3. $N$ is a symmetric endomorphism such that

\begin{equation} \label{eq8}
g(N\vec X,\vec Y)=\frac{1}{2}L_{\vec \xi}g(\vec X,\vec Y),
\end{equation}
where $\vec X,\vec Y,\vec Z\in\Gamma D$ are sections of the
distribution $D$, and $p:TX\to D$ is the projection.
\end{theorem}

{\bf Proof.} The first twoc onditions of the theorem uniquely
define the interior metric connection [7]. Alternating the second
covariant derivative we get
$$\nabla_{[e}\nabla_{a]}g_{bc}=2\omega_{ea}\partial_ng_{bc}-g_{dc}R^d_{eab}-g_{bd}R^d_{eac}.$$
Comparing the obtained equality with \eqref{eq8} we find the
implicit expression for the endomorphism $N$:
$$N^f_b=\frac{1}{4(n-1)}\omega^{ea}(R^f_{eab}+g_{bd}g^{cf}R^d_{eac}).$$
If $\partial_n g_{ab}-0$, then $N=0$. This proves the theorem.
$\Box$

We call the prolonged connection with the properties form theorem
2 the $N$-prolonged metric connection. We will use the notation
$\nabla^N=(\nabla,N)$ for the prolonged connection. In particular,
$\nabla^1=(\nabla, 0)$.

\subsection[2.3]{Special connections on manifolds with almost
compact metric structure} E.Cartan [8--10] was the first who
considered linear metric connection with a torsion instead of the
Levi-Civita connection. The most interesting among the metric
connections with torsion is the semi-symmetric connection
investigated by K.Yano in [11]. The quarter-symmetric connection
defined in 1975 Golab [12]. There is a big number of works
dedicated to metric and non-metric connections with torsion
defined on manifolds with almost contact structures. Here we fix
attention to the paper by Bejancu [13]. Bejancu defines the
connection $\nabla^B$ on a Sasaki manifold by the formula
$$\nabla^B_{\vec X}\vec Y=\tilde\nabla_{\vec X}\vec Y-\eta(\vec X)\tilde\nabla_{\vec Y}\vec
\xi-\eta(\vec Y)\tilde\nabla_{\vec X}\vec \xi+(\omega+c)(\vec
X,\vec Y)\vec\xi.$$ With respect to an adapted coordinates, the
non-zero components of the connection $\nabla^B$ are
$$\Gamma^{Ba}_{\,\,bc}=\Gamma^a_{bc}=\frac{1}{2}g^{ad}(\vec e_b g_{cd}-\vec e_c g_{bd}-\vec e_d
g_{bc}).$$ The constructed by Bejancu connection is in general not
metric in the common case of almost contact structure more general
then the Sasaki structure. Indeed, by the equality
$$\nabla^B_ng_{ab}=\partial g_{ab},$$
the metric Bejancu connection is equivalent to  an almost
K-contact  almost contact metric structure. Define on a manifold
with an almost contact metric structure the connection $\nabla^N$
by the equality $$\nabla^N_{\vec X}\vec Y=\nabla^B_{\vec X}\vec
Y+\eta(\vec X)N\vec Y,$$ where $N$ is the endomorphism from
theorem 2. Let us call the introduced connection the
$N$-connection. The non-zero components of this connection are
$$\Gamma^{Na}_{\,\,bc}=\Gamma^a_{bc}=\frac{1}{2}g^{ad}(\vec e_b g_{cd}-\vec e_c g_{bd}-\vec e_d
g_{bc}),quad \Gamma^{Na}_{\,\,nc}=N^a_c.$$ The torsion of the
$N$-connection is defined by the equality
$$S^N(\vec X,\vec Y)=2\omega(\vec X,\vec Y)\xi+\eta(\vec X)N\vec
Y-\eta(\vec Y)N\vec X.$$ The following theorem can be proved by
direct computations.

\begin{theorem} An $N$-connection is a metric
connection.\end{theorem}

\section[3]{$N$-prolonged connection as an almost contact metric
structure}

Consider on a manifold $X$ a contact metric structure
$(D,\varphi,\xi,\eta,g,X)$. We define on the distribution $D$ as
on a smooth manifold the almost contact metric structure $(\tilde
D,J,\vec u,\lambda=\eta\circ\pi_*,\tilde g, D)$ by setting
$$\tilde g(\vec\epsilon_a,\vec\epsilon_b)=\tilde g(\partial_{n+a},\partial_{n+b})=\tilde g(\vec e_a,\vec e_b),
\quad \tilde g(\vec\epsilon_a,\partial_{n+b})=\tilde
g(\vec\epsilon_a,\vec u)=\tilde g(\vec u,\partial_{n+b})=0,$$
 $$J(\vec\epsilon_a)=\partial_{n+a},\quad J(\partial_{n+a})=-\vec\epsilon_a,\quad J(\vec u)=0.$$
 The vector fields
 $$\vec\epsilon_a=\partial_a-\Gamma^n_a\partial_n-\Gamma^b_{ac}x^{n+c}\partial_{n+b},\quad\vec u=\partial_n-N^a_bx^{n+b}\partial_{n+a}$$ are defined here by the prolonged connection.
 Let $\tilde \omega=d\lambda$. It can be directly checked that the
 non-zero components of the form $\tilde \omega$ are given by the
 equality $\tilde\omega_{ab}=\omega_{ab}$. Consequently, ${\text
 rk}\tilde\omega=\frac{n-1}{2}$. This  implies that
 the constructed structure is not contact and, in particular, it
 is not a Sasaki structure.

\begin{theorem} The prolonged almost contact metric structure is
al most K-contact if and only if the initial structure is
K-contact.\end{theorem}

{\bf Proof.} Possibly non-zero components of the curvature tensor
with respect to adapted coordinates are of the form
\begin{align}\label{eq9} (L_{\vec u}\tilde g)_{ab}&=\partial_n
g_{ab},\\
\label{eq10} (L_{\vec u}\tilde g)_{n+a,n+b}&=\partial_n
g_{ab}-g_{ac}N^c_b-g_{cb}N^c_a,\\
\label{eq11} (L_{\vec u}\tilde
g)_{n+a,b}&=g_{ac}(P^c_{bd}-\nabla_bN^c_d)x^{n+d}.
\end{align}
In fact, the components from \eqref{eq10} are also zero, since
these are the components of the covariant derivative of the metric
tensor. The equality $\partial_ng_{ab}$ implies $N_d^c=0$ and
$P^c_{bd}=0$ (see \eqref{eq9} and \eqref{eq10}). This proofs the
theorem. $\Box$

Suppose now the=at the initial structure is K-contact ($N=0$),
then the following theorem takes a place.

\begin{theorem} The almost contact metric structure
$(\tilde D,J,\vec u,\lambda=\eta\circ\pi_*,\tilde g, D)$ is is
almost normal if and only if the distribution $D$ is a
distribution of zero curvature.\end{theorem}

{\bf Proof.} Let us rewrite the equality \eqref{eq1} in new
notation, $$N_{J}+2(d\tilde\eta \circ J)\otimes \vec{u}=0.$$ In
[4] it was shown that an almost contact structure is almost normal
if and only if $\tilde P\circ N_J=0$, where $\tilde P:TD\to\tilde
D$ is the projection.

Using the equalities \eqref{eq3}--\eqref{eq5} for the case of the
connection $\nabla^1$, we get the following expressions for the
Nijenhuis torsion of the operator $J$:
\begin{align*}
N_J\left({\vec{\varepsilon }}_a,{\vec{\varepsilon }}_b\right)&=-R^e_{abc}x^{n+c}{\partial }_{n+e},\\
N_J\left({\partial }_{n+a},{\partial }_{n+b}\right)&=2{\omega }_{ba}{\partial }_n+R^e_{abc}x^{n+c}{\partial }_{n+e},\\
N_J\left({\vec{\varepsilon }}_a,{\partial }_{n+b}\right)&=0,\\
N_J\left({\vec{\varepsilon }}_a,{\partial }_n\right)&=
N_J\left({\partial }_{n+a},{\partial
}_n\right)=-x^{n+c}P^b_{ac}{\partial }_{n+b}.\end{align*} Thus the
prolonged almost contact metric structure is almost normal if and
only if the Schouten curvature tensor is zero. $\Box$

\textbf{References}

1. S. Sasaki, {\it On the differential geometry of tangent bundles
of Riemannian manifolds.} Tohoku Math. J. 10 (1958), 338-354.

2. A. V. Bukusheva, S. V. Galaev, {\it Almost contact metric
structures defined by connection over distribution with admissible
Finslerian metric}. Izv. Saratov. Univ. Mat. Mekh. Inform., 12:3
(2012), 17–22

3. A. V. Bukusheva, S. V. Galaev, {\it    Connections over a
distribution and geodesic sprays.} Izv. Vyssh. Uchebn. Zaved.
Mat., 2013, no. 4,  10--18

4. S.V. Galaev, {\it The intrinsic geometry of almost contact
metric manifolds}, Izv. Saratov. Univ. Mat. Mekh. Inform., 12:1
(2012), 16--22.

5. S.V. Galaev, {\it Almost contact K\"ahlerian manifolds of
constant holomorphic sectional curvature.} Izv. Vyssh. Uchebn.
Zaved. Mat., 2014, no. 8,  42--52

6. S.V. Galaev,  A.V Gokhman, {\it About first integrals of
dynamical system with integrable linear connection}. Math. Mech.
Proc. Saratov univ. 2013. Vol. 15. P. 23--26.

7. V.V.~Wagner, {\it Geometry of  $(n - 1)$-dimensional
nonholonomic manifold
 in an  $n$-dimensional space},  Proc. Sem. on vect. and tens. analysis (Moscow
 Univ.).  5 (1941), 173--255.

8. Cartan E. Sur les vari?et?es `a connexion affine et la th?eorie
de la relativ?e. I // Ann. Sci. ?Ecole Norm. Sup.---1923.---Vol.
40.---P. 325---412.

9. Cartan E. Sur les vari?et?es `a connexion affine et la th?eorie
de la relativ?e g?en?eralis?ee. I // Ann. Sci. ?Ecole Norm.
Sup.---1924.---Vol. 41.---P. 1---25.

10. Cartan E. Sur les vari?et?es `a connexion affine et la
th?eorie de la relativ?e g?en?eralis?ee. II // Ann. Sci. ? Ecole
Norm. Sup.---1925.---Vol. 42.---P. 17---88.

11. Yano K. On semi-symmetric metric connections, Revue Roumaine
de Math. Pures et Appliques 15 (1970), 1579-1586.

12. Golab S. On semi-symmetric and quarter-symmetric linear
connections, Tensor N.S., 29(1975), 249-254.

13. Bejancu A. Kahler contact distributions. Journal of Geometry
and Physics 60, 1958--1967 (2010).

\vskip1cm

Saratov State University, Department of Mechanic and Mathematics,  Chair of Geometry\\
Astrakhanskaya 83, 410012 Saratov, Russia

sgalaev@mail.ru

\end{document}